\documentclass[12pt]{amsart}
\usepackage{amssymb,amsmath,mathrsfs}

\usepackage[T1]{fontenc}

\newcommand{\N}{{\mathbb N}}

\newcommand{\R}{{\mathbb R}}

\newcommand{\eps}{\varepsilon}

\numberwithin{equation}{section}
\newtheorem{theorem}{Theorem}[section]
\newtheorem{proposition}[theorem]{Proposition}
\newtheorem{lemma}[theorem]{Lemma}

\theoremstyle{definition}
\newtheorem{remark}[theorem]{Remark}

\newcommand{\brm}{\begin{remark}\rm}
\newcommand{\erm}{\end{remark}}
\newcommand{\brms}{\begin{remark}\rm}
\newcommand{\erms}{\end{remark}}
\newcommand{\bte}{\begin{theorem}}
\newcommand{\ete}{\end{theorem}}
\newcommand{\bpr}{\begin{proposition}}
\newcommand{\epr}{\end{proposition}}
\newcommand{\ble}{\begin{lemma}}
\newcommand{\ele}{\end{lemma}}
\newcommand{\beq}{\begin{equation}}
\newcommand{\eeq}{\end{equation}}
\newcommand{\bdm}{\begin{displaymath}}
\newcommand{\edm}{\end{displaymath}}
\numberwithin{equation}{section}

\newcommand{\bos}{\begin{remark}\rm}
\newcommand{\eos}{\end{remark}}

\newcommand{\ben}{\begin{enumerate}}
\newcommand{\een}{\end{enumerate}}

\newcommand{\be}{\begin{equation}}
\newcommand{\ee}{\end{equation}}

\title[On a general class of Schr\"odinger systems II]{On Schr\"odinger systems with \\ local and nonlocal nonlinearities - Part2}

\author[H.\ Hajaiej]{Hichem Hajaiej}

\thanks{Ipeit (Institut preparatoire aux etudes d'ingenieur de Tunis)
 2, Rue Jawaher Lel Nahru -
1089 Montfleury - Tunis,
 Tunisie.
E-mail: {\em hichem.hajaiej@gmail.com}}

\address{Ipeit (Institut preparatoire aux etudes d'ingenieur de Tunis)
\newline\indent
 2, Rue Jawaher Lel Nahru -
1089 Montfleury - tunis
 Tunisie}
\email{hichem.hajaiej@gmail.com}

\thanks{The author was partially supported by the
 the Tunisian ARUB project : {\em Analyse
Math\'ematique et Applications 04/UR/15-02. }}

\begin{document}

\subjclass[2000]{35J40; 58E05}

\keywords{Nonlinear Schr\"odinger equations, soliton dynamics in external potential,
ground states, semi-classical limit.}

\begin{abstract}
In this second  part, we establish the  existence of  special solutions of the nonlinear Schr\"odinger system  studied  in the  first part when the  diamagnetic  field  is nul. 
We also prove  some symmetry properties  of  these  ground  states  solutions.
\end{abstract}
\maketitle

\medskip
\begin{center}
\begin{minipage}{11cm}
\footnotesize
\tableofcontents
\end{minipage}
\end{center}
\medskip

\section{Study of ground state solutions}

\subsection{Introduction and historical remarks}
In this section, we shall study the existence and symmetry of ground states for the following
$m \times m$ nonlinear Schr\"odinger system without magnetic field, in presence of local an nonlocal
nonlinearities
\begin{equation}\label{Nomagnetic}
\begin{cases}
- \Delta \Phi_j + (\lambda-V(|x|)) \Phi_j - g_j(|x|, |\Phi_1|^2,\dots, |\Phi_m|^2) \Phi_j - \sum\limits_{i=1}^m W_{ij} * h(|\Phi_i|)  \frac{h'(|\Phi_j|)}{|\Phi_j|} \Phi_j=0 &\\
\noalign{\vskip4pt}
\,\,\text{for $1 \leq j \leq m$}.
\end{cases}
 \end{equation}
For every $\Phi=(\Phi_1, \dots, \Phi_m) \in \mathcal{H}^1(\R^N)$, we define the energy functional
\begin{align*}
{\mathcal E}(\Phi)& =\frac12 \sum_{j=1}^ m \int \left| \nabla \Phi_j  \right|^2  dx
- \frac12 \int V(|x|) |\Phi|^2 \, dx- \int G(|x|, |\Phi_1|^2,\dots,|\Phi_m|^2)  dx   \\ \nonumber
& -  \frac12 \sum_{i, j=1}^{m} \iint W_{ij}(|x-y|)  h(|\Phi_{i}(x)|)  h(|\Phi_{j}(y)|)  dx dy.
\end{align*}
We are interested to solve the following minimization problem
\begin{align}
    \label{minimizatProbl}
I_c= \inf_{\Phi \in {\mathcal S}_c} {\mathcal E}(\Phi),\qquad
{\mathcal S}_c=\Bigl\{  \Phi \in \mathcal{H}^1(\R^N): \  \sum\limits_{j=1}^ m \int | \Phi_j|^2=c \Bigr\},
\end{align}
where $c > 0$ is a fixed number.

\section{Main assumptions}

\subsection{Assumptions on local nonlinearities}

We assume that the following conditions hold
\vskip2pt
\noindent
$(V0)$ $V: \R^N \rightarrow \R^+$ satisfies
$$
V(|x|) \geq V(|y|), \quad\text{for all $x, y \in \R^N$ with $ |x| \leq |y|$ }.
$$
Moreover,
$$
V(|x|) \rightarrow 0, \quad\hbox{as $|x| \rightarrow \infty$.}
$$

\vskip2pt
\noindent
$(G0)$  $G: (0, \infty) \times \R^m \to \R$ is a super-modular function, namely
\begin{align}\label{coop1}
G(r, y+h e_i+k e_j) +G(r,y) &\geq G(r, y+h e_i) +G(r,y+ke_j) \\
\label{coop2}
G(r_1, y+h e_i) +G(r_0,y) &\leq G(r_1, y) +G(r_0,y+he_i)
\end{align}
for $i \neq j$, $h,k >0$, $y=(y_1,\dots,y_m)$ and $\{e_i\}$ is the
standard basis in $\R^m$, $r >0$ and $0 < r_0 <r_1$.
\vskip2pt
\noindent
$(G1)$  There exists $K>0$ such that, for all $r>0$ and $s_1,\dots,s_m  \geq 0$, we have
$$
0 \leq G(r, s_1,\dots, s_m) \leq K \Big(\sum_{j=1}^{m} s_j + \sum_{j=1}^{m} s_j^{\frac{\ell_j +2}{2}}   \Big), \qquad  0 < \ell_j <  \frac4N.
$$
\vskip2pt
\noindent
$(G2)$ for all $\varepsilon >0$, there exist $R_0 >0$ and $S_0 >0$ such that
$G(r, s_1,\dots, s_m) \leq \varepsilon \sum_{j=1}^m s_j$, for all $r > R_0$ and $s_1,\dots,s_m < S_0$;
\vskip2pt
\noindent
$(G3)$ For any $r > 0$, $s_1, \dots, s_m$ and $t > 1$,
$$
G(r, t s_1,\dots, t s_m) \geq t G(r, s_1,\dots, s_m).
$$
\vskip2pt
\noindent
$(G4)$  There exist $B, \gamma, R_2, S_2  > 0$  such that
$$
G(r, s_1, 0, \dots, 0) \geq B s_1^{\gamma}, \quad\hbox{for any $r > R_2$, $0 \leq s_1 \leq S_2$},
$$
where $1 \leq \gamma < 1+ \frac 2N$.

\subsection{Assumptions on the nonlocal nonlinearities}

\noindent We need the following  assumptions
\smallskip

\noindent $(h0)$ $h: \R^+ \to \R^+$ is continuous, non-decreasing, $h(0)=0$ and there exists $M > 0$ such that
$$
h(s) \leq M s^{\mu} \quad \hbox{where \quad $2 \leq \mu < 2 - \frac{1}{q}+ \frac{2}{N} $};
$$

\medskip

\noindent $(h1)$ $h(t s) \geq  t h(s)$, \,\, for all $t > 1$ and $s \geq 0$.
\medskip

\noindent  $(h2)$  There exist   $A,  S_1 > 0 $ and $\beta \geq \mu$   such that $h(s) \geq A s^{\beta}$, for any $ 0 \leq  s  \leq S_1$.
\medskip

\noindent $(W1)$ There exist  $\Gamma, C, t_1  > 0$ such that
$$
W_{1 1}\left(\frac rt \right) \geq C \frac{t^{\Gamma}}{r^{\Gamma}}, \quad\hbox{for any $r \geq 0$, $0 \leq t \leq t_1$}
$$
where  $2N-N \beta -\Gamma +2 > 0$.

\section{Sign of the Lagrange multiplier}

We have the following

\begin{proposition}
Let $c>0$ and assume that the minimization problem~\eqref{minimizatProbl} admits a solution $\hat\Phi\in{\mathcal S}_c$
with negative energy, namely
$$
{\mathcal E}(\hat\Phi)=I_c<0.
$$
Assume furthermore that the function
\begin{equation*}
N(\Phi)=\int G(|x|, |\Phi_1|^2,\dots,|\Phi_m|^2) dx
+\frac12 \sum_{i, j=1}^{m} \iint W_{ij}(|x-y|)  h(|\Phi_{i}(x)|)  h(|\Phi_{j}(y)|)dxdy
\end{equation*}
satisfies over $\hat\Phi$ the condition
\begin{equation}
    \label{segnoNN}
N'(\hat\Phi_1,\dots,\hat\Phi_m)(\hat\Phi_1,\dots,\hat\Phi_m)-2N(\hat\Phi_1,\dots,\hat\Phi_m)\geq 0.
\end{equation}
Let $\lambda_c$ denote the Lagrange multiplier associated with $\hat\Phi$. Then $\lambda_c<0$.
\end{proposition}

\begin{proof}
Of course, we have ${\mathcal E}'(\hat\Phi)=\lambda_c \hat \Phi$, so that
$$
{\mathcal E}'(\hat\Phi)(\hat \Phi)=\lambda_c (\hat \Phi,\hat\Phi)_{\mathcal{L}^2}=\lambda_c \|\hat\Phi\|_{\mathcal{L}^2}^2= c\lambda_c.
$$
Then, we have
\begin{equation*}
c\lambda_c-2 I_c  ={\mathcal E}'(\hat\Phi)(\hat \Phi)-2{\mathcal E}(\hat\Phi)
= -N'(\hat\Phi)(\hat\Phi)+2N(\hat\Phi)=\tau,
\end{equation*}
namely $\lambda_c=\frac{2I_c}{c}+\frac{\tau}{c}<0$, as $\tau \leq 0$ and $I_c<0$ by assumption.
This proves the assertion.
\end{proof}

\begin{remark}
Assume that the function $\R^m\ni s\mapsto G(r,s)\in\R^+$ is homogeneous of degree $\varrho\geq 1$
and $W_{ij}(x)=0$ for all $i,j=1,\dots,m$ and $x\in\R^N$. Then condition~\eqref{segnoNN} is satisfied. In fact, taking into
account that $\nabla G(s)\cdot s=dG(s)(s)=\varrho G(s)$, it follows that
\begin{align*}
N'(\hat\Phi)(\hat\Phi)-2N(\hat\Phi)& =2\int \sum_{j=1}^m D_{s_j}G(|x|, |\hat\Phi_1|^2,\dots,|\hat\Phi_m|^2)||\hat\Phi_j|^2dx \\
&-2\int G(|x|, |\hat\Phi_1|^2,\dots,|\hat\Phi_m|^2)dx \\
& = 2(\varrho-1)\int G(|x|, |\hat\Phi_1|^2,\dots,|\hat\Phi_m|^2)dx \geq 0,
\end{align*}
which proves the desired claim. The homogeneity of $G$ is often fulfilled in the applications. Think,
instance, to the literature of weakly coupled nonlinear Schr\"odinger systems.
\end{remark}

\begin{remark}
Assume that the function $s\mapsto h(s)$ is homogeneous of degree $\mu\geq 2$ and that $G=0$.
Then condition~\eqref{segnoNN} is satisfied. In fact, taking into
account that $h'(s)s=\mu h(s)$, by direct computation, exchanging $i$ and $j$
and $x$ with $y$, it follows that
\begin{align*}
N'(\hat\Phi)(\hat\Phi)-2N(\hat\Phi)& =
\frac{1}{2}\sum_{i, j=1}^{m} \iint W_{ij}(|x-y|)  h(|\hat\Phi_{i}(x)|)  h'(|\hat\Phi_{j}(y)|)|\hat\Phi_j(y)|dxdy \\
& +\frac{1}{2}\sum_{i, j=1}^{m} \iint W_{ij}(|x-y|)  h(|\hat\Phi_{j}(y)|)  h'(|\hat\Phi_{i}(x)|)|\hat\Phi_i(x)|dxdy \\
& -\sum_{i, j=1}^{m} \iint W_{ij}(|x-y|)  h(|\hat\Phi_{i}(x)|)  h(|\hat\Phi_{j}(y)|)dxdy \\
& =\sum_{i, j=1}^{m} \iint W_{ij}(|x-y|)  h(|\hat\Phi_{i}(y)|)  h'(|\hat\Phi_{j}(x)|)|\hat\Phi_j(x)|dxdy \\
& -\sum_{i, j=1}^{m} \iint W_{ij}(|x-y|)  h(|\hat\Phi_{i}(x)|)  h(|\hat\Phi_{j}(y)|)dxdy \\
& =(\mu-1)\sum_{i, j=1}^{m} \iint W_{ij}(|x-y|)  h(|\hat\Phi_{i}(x)|)  h(|\hat\Phi_{j}(y)|)dxdy\geq 0,
\end{align*}
which proves the claim. The homogeneity of $h$ is often fulfilled in the applications. Think for
instance to the literature of the {\em Pekar-Choquard} equation with $h(s)=|s|^\mu$, being the classical
formulation in the particular case $\mu=2$.
\end{remark}

\section{Existence and symmetry of solutions}

We have the following

\begin{proposition}
    \label{welpMP}
Assume conditions (V0), (G1),  (h0) hold. Then, for all $c>0$, problem~\eqref{minimizatProbl} is well-posed,
that is $I_c>-\infty$.
\end{proposition}

\begin{proof}
    Let $\Phi \in {\mathcal S}_c$. In the following, we shall denote by $C$ a generic positive constant, possibly depending on $c$,
that can change from line to line. From assumption (G1), we have
\begin{equation}
    \label{localineqq}
\int G(|x|, |\Phi_1|^2,\dots,|\Phi_m|^2)dx \leq  C + C\sum_{j=1}^{m} \| \Phi_j\|^{\ell_j+2}_{\ell_j+2}.
\end{equation}
From the Gagliardo-Nirenberg inequality, and since $\|\Phi_j \|_{L^2}\leq \sqrt{c}$, we have
$$
\| \Phi_j   \|_{\ell_j +2}^{\ell_j+2} \leq C\| \Phi_j \|^{(1- \sigma_j)(\ell_j+2)}_{L^2}
\|\nabla  \Phi_j \|_{L^2}^{\sigma_j(\ell_j+2)}\leq C\|\nabla  \Phi_j \|_{L^2}^{\sigma_j(\ell_j+2)}, \qquad
\sigma_j= \frac{N \ell_j}{2(\ell_j +2)},
 $$
for $j=1,\dots,m$. Notice that, by assumption, we have
$$
\sigma_j(\ell_j+2)=\frac{N\ell_j}{2}<2,\quad \text{for $j=1,\dots,m$}.
$$
Then, by means of Young inequality, for all $\eps>0$ there exists $K_1(\eps)>0$ such that
\begin{align}\label{youngIneq}
\| \Phi_j   \|_{\ell_j +2}^{\ell_j +2}  \leq  K_1(\eps)+\eps\| \nabla  \Phi_j \|_{L^2}^2.
\end{align}
In turn, inequality~\eqref{localineqq} yields
\begin{equation}
    \label{loc1fin}
    \int G(|x|, |\Phi_1|^2,\dots,|\Phi_m|^2)dx \leq K_1(\eps) + \eps\sum_{j=1}^{m} \| \nabla  \Phi_j \|_{L^2}^2,
\end{equation}
for some positive constant $K_1(\eps)$. Dealing with the nonlocal nonlinearities,
from assumption $(h)$, by the Hardy-Littlewood inequality combined with
the Gagliardo-Nirenberg inequality, for any $i, j= 1, \dots, m$, since
$\max\{\| W_{ij} \|_{L^q_w}:i,j=1,\dots,m\}<\infty$, setting
$$
\hat q=\frac{2q}{2q-1},\qquad
\gamma= \frac{N}{2} \left(\frac{\hat q \mu-2}{\hat q \mu}\right),
$$
for every $\eps>0$ there exists $K_2(\eps)>0$ such that
\begin{align}
    \label{nonlocccfin}
&\frac{1}{2}\sum_{i,j=1}^m \iint W_{ij}(|x-y|)  h(|\Phi_{i}(x)|)  h(|\Phi_{j}(y)|) dx dy
\leq C\sum_{i,j=1}^m \| W_{ij} \|_{L^q_w} \| \Phi_i ^{\mu}\| _{L^{\hat q}} \| \Phi_j^{\mu}\|_{L^{\hat q}}   \\
& \leq C\sum_{i,j=1}^m \| \Phi_i\|_{L^{\hat q\mu}}^\mu \| \Phi_j\|_{L^{\hat q\mu}}^\mu \leq
C\sum_{i,j=1}^m \| \Phi_i\|_{L^2}^{(1-\gamma)\mu} \| \nabla \Phi_i\|_{L^2}^{\gamma\mu}
\| \Phi_j\|_{L^2}^{(1-\gamma)\mu} \| \nabla \Phi_j\|_{L^2}^{\gamma\mu} \notag \\
&\leq C\sum_{i,j=1}^m \| \nabla \Phi_i\|_{L^2}^{\gamma\mu}\| \nabla \Phi_j\|_{L^2}^{\gamma\mu} \leq
C\sum_{i=1}^m \| \nabla \Phi_i\|_{L^2}^{2\gamma\mu}\leq K_2(\eps)+\eps \sum_{i=1}^m \| \nabla \Phi_i\|_{L^2}^{2}, \notag
\end{align}
where in the last two inequalities we used the Young inequality. In particular, the last one
was possible since, by our assumptions on $\mu$ in $(h0)$, we have
$$
2\gamma \mu=N\left(\frac{\hat q \mu-2}{\hat q}\right)=
N\left(\frac{q \mu-2q+1}{q}\right)<2.
$$
Then, fixed $\eps\in (0,1/4)$, by combining~\eqref{loc1fin} and~\eqref{nonlocccfin}, by the definition of ${\mathcal E}$ and denoted by $\rho= V(0) > 0$, we have
\begin{align}\label{Functbded}
{\mathcal E}(\Phi)& \geq \frac12\sum_{j=1}^ m \| \nabla \Phi_j  \|^2_{L^2}
- \frac{\rho}{2} \sum_{j=1}^ m\|\Phi_j\|^2_{L^2}- \int G(|x|, |\Phi_1|^2,\dots,|\Phi_m|^2)  dx   \\ \nonumber
& -  \frac12 \sum_{i, j=1}^{m} \iint W_{ij}(|x-y|)  h(|\Phi_{i}(x)|)  h(|\Phi_{j}(y)|)  dx dy \\
&\geq \Big(\frac12-2\eps\Big) \sum_{j=1}^ m \| \nabla \Phi_j  \|^2_{L^2}
- \frac{\rho c}{2} -K_1(\eps)-K_2(\eps)\geq - \frac{\rho c}{2}-K_1(\eps)-K_2(\eps).
\end{align}
for all $\Phi\in {\mathcal S}_c$, yielding the desired conclusion.
\end{proof}

The next proposition shows that, even in the limiting cases with respect to the growths
of the local and nonlocal nonlinearities the minimization problem is well posed, provided
that the infimum is taken over a sphere of sufficiently small radius $c$.

\begin{proposition}
Assume conditions (V0), (G1), (h0) hold and that
$$
\text{either $\ell_{j_0}=\frac{4}{N}$ for some $j_0=1,\dots,m$ or $\mu=2-\frac{1}{q}+\frac{2}{N}$}.
$$
Then $I_c>-\infty$ for every $c>0$ sufficiently small.
\end{proposition}
\begin{proof}
Let $c>0$ and take $\Phi \in {\mathcal S}_c$. In the following, we shall denote
by $C$ a generic positive constant which can change from line to line and which
is independent of $c$. In fact, differently from the proof of Proposition~\ref{welpMP},
here we need to put $c$ into evidence in the estimates in order to show that
problem~\eqref{minimizatProbl} is well posed, for all $c$ sufficiently small.
Assume that there exists $1 \leq j_0 \leq m$ such that $ \ell_{j_0} =\frac4N$ (and that
$\ell_j<4/N$ for all $j\neq j_0$). Recall that
$\|\Phi_{j_0} \|_{L^2}\leq \sqrt{c}$. From (G1), the Gagliardo-Nirenberg
inequality and~\eqref{youngIneq} (holding, indeed, when $\ell_j<4/N$), we have
\begin{align*}
\label{clocalineqq}
\int G(|x|, |\Phi_1|^2,\dots,|\Phi_m|^2)dx & \leq  C +  C \| \Phi_{j_0}   \|_{\ell_{j_0} +2}^{\ell_{j_0}+2}     +C \sum_{j \neq j_0}^{m} \| \Phi_j\|^{\ell_j+2}_{\ell_j+2} \\
& \leq  C  +  C  \| \Phi_{j_0} \|^{\frac{4}{N}}_{L^2} \|\nabla  \Phi_{j_0} \|_{L^2}^2 +   K_1(\eps) + \eps \sum\limits_{j \neq j_0}^m
\|\nabla  \Phi_{j} \|_{L^2}^2  \\
& \leq   K_1(\eps)+   C c^{\frac{2}{N}} \|\nabla  \Phi_{j_0} \|_{L^2}^2  + \eps \sum\limits_{j \neq j_0}^m
\|\nabla  \Phi_{j} \|_{L^2}^2  \\
& \leq   K_1(\eps)+   \max\{C c^{\frac{2}{N}},\eps\}\sum\limits_{j=1}^m
\|\nabla  \Phi_{j} \|_{L^2}^2
\end{align*}
for some positive constant $K_1(\eps)$ depending on $\eps$.
Concerning the nonlocal non\-li\-nea\-ri\-ties, we observe that,  if $\mu < 2-1/q+2/N$, we are in the case
of the proof of Proposition~\ref{welpMP} and we have inequality~\eqref{nonlocccfin}.
If, instead, we are in the limiting case $\mu = 2-1/q+2/N$, for  $\hat q=\frac{2q}{2q-1}$ it holds
$$
\gamma= \frac{1}{\mu}= \frac{N q}{2Nq-N+2q}.
$$
In turn, by Hardy-Littlewood and Gagliardo-Nirenberg inequalities, we have
\begin{align*}
&  \frac{1}{2}\sum_{i,j=1}^m \iint W_{ij}(|x-y|)  h(|\Phi_{i}(x)|)  h(|\Phi_{j}(y)|) dx dy  \\
           &   \leq
C\sum_{i,j=1}^m \| \Phi_i\|_{L^2}^{(1-\gamma)\mu} \| \nabla \Phi_i\|_{L^2}
\| \Phi_j\|_{L^2}^{(1-\gamma)\mu} \| \nabla \Phi_j\|_{L^2} \notag \\
& \leq  C c^{(1-\gamma) \mu} \sum_{i,j=1}^m \| \nabla \Phi_i\|_{L^2}\| \nabla \Phi_j\|_{L^2}
 \leq
C c^{(1-\gamma) \mu} \sum_{i=1}^m \| \nabla \Phi_i\|_{L^2}^{2}. \notag
\end{align*}
In any case, by~\eqref{nonlocccfin} and the above inequality, we can always write
$$
\frac{1}{2}\sum_{i,j=1}^m \iint W_{ij}(|x-y|)  h(|\Phi_{i}(x)|)  h(|\Phi_{j}(y)|) dx dy
\leq \max\{C c^{(1-\gamma) \mu},\eps\}\sum_{i=1}^m \| \nabla \Phi_i\|_{L^2}^{2}+K_2(\eps).
$$
Then, by the definition of ${\mathcal E}$ and previous inequalities, denoted by $\rho= V(0) > 0$, we have
\begin{align*}
{\mathcal E}(\Phi)& \geq    \Big(\frac12 - \max\{C c^{\frac{2}{N}},\eps\} - \max\{C c^{(1-\gamma) \mu},\eps\}  \Big) \sum_{j=1}^m \| \nabla \Phi_j\|_{L^2}^{2}  - \frac{\rho c}{2}-K_1(\eps)-K_2(\eps),
\end{align*}
for all $\Phi \in {\mathcal S}_c$. By choosing $\eps>0$ and $c>0$ so small that
$$
\frac12 - \max\{C c^{\frac{2}{N}},\eps\} - \max\{C c^{(1-\gamma) \mu},\eps\}>0
$$
it holds ${\mathcal E}(\Phi)\geq - \frac{\rho c}{2}-K_1(\eps)-K_2(\eps)$ and the assertion follows,
namely there exists $c_0>0$ such that the minimization problem is well posed for all $c\in (0,c_0)$.
\end{proof}

The next proposition says that, at least under suitable assumptions, which include some
classical situations, such as $h(s)=s^\mu$, $W_{ij}(x)=|x|^{-\alpha}$ and
$$
G(|x|,s_1,\dots,s_m)=\frac{1}{\ell+2}\sum_{i,j=1}^m |s_i|^{(\ell+2)/2}+2|s_i|^{(\ell+2)/4}|s_j|^{(\ell+2)/4},
$$
the upper bounds on $\ell_j$ and $\mu$ are optimal
for the minimization problem to be well posed.

\begin{proposition}
Assume ($V0$) and that either there exists
a function $H:\R^m_+\to\R$, homogeneous of degree $\frac{\ell+2}{2}$ with $\ell>4/N$, such that
$$
G(|x|,s_1,\dots,s_m)\geq H(s_1,\dots,s_m),\quad\text{for all $(s_1,\dots,s_m)\in\R^m_+$}
$$
or there exist two constants $\gamma_1,\gamma_2>0$ such that, for some $1 \leq i_0, j_0 \leq m$,
$$
\text{$W_{ i_0 j_0}(x)\geq \gamma_1 |x|^{-\alpha}$ and $h(s)\geq\gamma_2 s^\mu$ for all $x\in\R^N$ and $s\in\R^+$,
with $\mu>2-\frac{\alpha}{N}+\frac{2}{N}$}.
$$
Then $I_c=-\infty$ for every $c>0$.
\end{proposition}
\begin{proof}
We consider the case when both the situations indicated in the statement occur, the proof being similar
in the other cases. Let $c>0$ and consider a fixed function $\Phi_0$ in ${\mathcal S}_c$.  For all $t>0$, we
define the function $\Phi_t:\R^N\to\R^m$ by setting $\Phi_t^j(x)=t^{N/2}\Phi_0^j(tx)$
for all $x\in\R^N$ and $j=1,\dots,m$. It follows that $\Phi_t\in {\mathcal S}_c$ for all $t>0$, so that,
by definition of $I_c$, it holds for all $t>0$ large
\begin{align*}
I_c\leq {\mathcal E}(\Phi_t)& \leq  \frac12\sum_{j=1}^ m \| \nabla \Phi_t^j  \|^2_{L^2}
- \int G(|x|, |\Phi_t^1|^2,\dots,|\Phi_t^m|^2)  dx   \\ \nonumber
& -  \frac12 \sum_{i, j=1}^{m} \iint W_{ij}(|x-y|)  h(|\Phi_{t}^i (x)|)  h(|\Phi_{t}^j (y)|)  dx dy  \\
& \leq  \frac{t^2}{2}\sum_{j=1}^ m \| \nabla \Phi_0^j  \|^2_{L^2}
       - t^{\frac{N \ell}{2}} \int H(|\Phi_0^1|^2,\dots,|\Phi_0^m|^2)  dx \\
&  - \frac{\gamma_1 \gamma_2^2}{2} t^{\alpha + N \mu -2N}  \iint |x-y|^{-\alpha} |\Phi_0^{i_0} (x) |^{\mu} |\Phi_0^{j_0} (y) |^{\mu} dx dy \\
& \leq C_1 t^2  - C_2 t^{\frac{N \ell}{2}}  - C_3 t^{\alpha + N \mu -2N} + C_4 \\
\noalign{\vskip3pt}
& \leq C_1 t^2 - C_5 t^{\min\{ \frac{N \ell}{2},  \alpha + N \mu -2N\}}+C_4.
\end{align*}
By assumptions $\min\{ \frac{N \ell}{2},  \alpha + N \mu -2N   \} > 2$ and the
assertion follows by letting $t\to\infty$.
\end{proof}

\begin{proposition}\label{SchwarzLimit}
Assume conditions ($V0$), ($W$), ($h0$), ($G0$), ($G1$) and ($G2$) hold.
Then, for every $c > 0$, problem~\eqref{minimizatProbl} admits a minimization
sequence $(\Phi_n)$ having a Schwarz symmetric weak limit $\Phi_0$ such that
$\mathcal E (\Phi_0) \leq I_c$ .
\end{proposition}

\begin{proof}
Let $\Phi_n \in \mathcal{H}^1$ be a minimizing sequence for~\eqref{minimizatProbl}.
Since $\|\nabla |\Phi_{n, j} | \|_{L^2}= \|   \nabla \Phi_{n, j}  \|_{L^2}$, we have that
$\mathcal E (|\Phi_n|)   \leq \mathcal E (\Phi_n) $ so that $|\Phi_n|$ is a minimizing sequence too.
In turn, without loss of generality, we may assume that the minimizing sequence is positive.
Denoted by $\Phi_n^*$ the sequence of the Schwarz symmetrizations of $\Phi_n$,  we claim that
$\mathcal E ( \Phi^*_n) \leq \mathcal E (|\Phi_n|)$ so that $ \Phi_n^* $
is also a minimizing sequence for~\eqref{minimizatProbl}.  In order to prove it, we take
advantage of the following symmetrization inequalities. By~\cite{HardyLittlePolya}, for every $j=1, \dots, m$,
 $$
  \|  \nabla  \Phi_{n, j}^*  \|_{L^2}^2 \leq \|  \nabla \Phi_{n, j}  \|_{L^2}^2
 $$
 $$
  \|   \Phi_{n, j}^*  \|_{L^2}^2=  \|   \Phi_{n, j}  \|_{L^2}^2.
 $$
From the last equality, it follows  that, if $\Phi_n \in {\mathcal S}_c$, then also $\Phi^*_n \in {\mathcal S}_c$.
Moreover, in view of assumption ($V0$),  we have that
$$
\int V(|x|)\Phi_{n, j}^2 \leq  \int V(|x|) (\Phi_{n, j}^*)^2.
$$
Furthermore, in view of the super-modularity assumption ($G0$), we have
    $$
    \int G(|x|, \Phi_{n, 1}^2,\dots,\Phi_{n, m}^2)dx \leq \int G(|x|, (\Phi_{n, 1}^*)^2,\dots,(\Phi_{n, m}^*)^2) dx
    $$
    and, by assumptions ($W$) and ($h0$), it follows
    $$
    \iint W_{ij}(|x-y|)  h(\Phi_{n, i}(x))  h(\Phi_{n, j}(y)) \, dx dy \leq
    \iint W_{ij}(|x-y|)  h(\Phi_{n, i}^*(x))  h(\Phi_{n, j}^*(y)) \, dx dy ,
    $$
    for every  any $i,j=1,\dots,m$.
We shall denote by $\widetilde{\Phi}_n=  \Phi_n^*$ a Schwarz symmetric minimizing sequence for~\eqref{minimizatProbl}.
Observe that $\widetilde{\Phi}_n$ is bounded in  $\mathcal{H}^1$. Indeed, if this was not the case, from the
following inequality (see inequality~\eqref{Functbded} in Proposition~\ref{welpMP}), as $n\to\infty$, denoted by $\rho= V(0) > 0$,
\begin{align*}
I_c + o(1) = {\mathcal E}(\widetilde{\Phi}_n) \geq
 \Big(\frac12-2\eps\Big) \sum_{j=1}^ m \| \nabla \widetilde{\Phi}_{n, j}  \|^2_{L^2}
- \frac{\rho c}{2} -K_1(\eps)-K_2(\eps)
\end{align*}
for $\eps \in (0, \frac{1}{4})$, we would immediately get a contradiction.
Hence, up to a subsequence, there exists $\Phi_0 \in  \mathcal{H}^1$ such that
 $\widetilde{\Phi}_n$ converges to $\Phi_0 $ weakly in $\mathcal{H}^1$,
locally strongly in $\mathcal{L}^s$ for  $s < 2^*$ and almost everywhere in $\R^N$.
 We will prove that
 \begin{equation}
    \label{semicontFunct}
 \mathcal E (\Phi_0) \leq \liminf_{n \rightarrow \infty} \mathcal E (\widetilde{\Phi}_n).
 \end{equation}
 For all $ j=1, \dots, m$, we know that
\begin{align}\label{gradSemicontInf}
\int  |\nabla \Phi_{0, j}|^2    \leq \liminf_{n \rightarrow \infty} \int  |\nabla \widetilde {\Phi}_{n, j}|^2.
\end{align}
Now, let us prove that, for every $i= 1, \dots, m$,
\begin{align}\label{convergV}
\lim_{n \rightarrow \infty}   \int V(|x|)  \widetilde{\Phi}_{n, j}^2=  \int V(|x|) \Phi_{0, j}^2,
\end{align}
\begin{align}\label{convergG}
\lim_{n \rightarrow \infty}  \int G(|x|, \widetilde{\Phi}_{n, 1}^2,\dots, \widetilde{\Phi}_{n, m}^2)=
\int G(|x|, \Phi_{0, 1}^2,\dots,\Phi_{0, m}^2),
\end{align}
and for all $i, j= 1, \dots, m$,
\begin{align}\label{convergW}
\lim_{n \rightarrow \infty} \iint W_{ij}(|x-y|)  h(\widetilde{\Phi}_{n, i}(x))  h(\widetilde{\Phi}_{n, j}(y))=
    \iint W_{ij}(|x-y|)  h(\Phi_{0, i}(x))  h(\Phi_{0, j}(y)).
\end{align}
First, we prove~\eqref{convergV}. Fixed $R > 0$,  denote by $B(R)$ the ball of radius $R$
centered at the origin. Since $\widetilde{\Phi}_{n, j} (x)\to\Phi_{0, j}(x) $ for a.e.~$x \in B(R)$
and there exists a function $b_j \in L^2(B(R))$ such that $\widetilde{\Phi}_{n, j} (x) \leq
b_j (x) $ for a.e.~$x \in B(R)$,  by the monotonicity assumption on $V$ in ($V0$), we have
\begin{align}\label{VBR}
\lim_{n \rightarrow \infty}   \int_{B(R)} V(|x|)  \widetilde{\Phi}_{n, j}^2=  \int_{B(R)} V(|x|) | \Phi_{0, j}|^2,
\end{align}
by dominated convergence. Now, fix $\eps > 0$ and $j=1,\dots,m$. Since $V(|x|)\to 0$ as $|x|\to\infty$ by assumption ($V0$),
there exists  $R(\eps)>0$ such that, for all $|x|>R(\eps)$ and for every $n \in \N$
$$
  \int_{B_c(R(\eps))} V(|x|)  \widetilde{\Phi}_{n, j}^2 \leq \eps     \int_{B_c(R(\eps))}  \widetilde{\Phi}_{n, j}^2 \leq \eps c.
$$
Furthermore, in a similar fashion, we have that
$$
  \int_{B_c(R(\eps))} V(|x|)  \widetilde{\Phi}_{0, j}^2(x)  \leq \eps c.
$$
By means of~\eqref{VBR}, choosing $R=R(\eps)$, there exists $\nu_{\eps} \in \N$
such that  for every  $n \geq \nu_{\eps} $
$$
\left|  \int_{B(R(\eps))} V(|x|)  \widetilde{\Phi}_{n, j}^2-  \int_{B(R(\eps))} V(|x|) \Phi_{0, j}^2   \right| < \eps.
$$
Thus, by combining the above inequalities, \eqref{convergV} follows. Now,
we show~\eqref{convergG}. Fixed $R  > 0$, it holds
 \begin{align}\label{GBR}
\lim_{n \rightarrow \infty}  \int _{ B(R)} G(|x|, \widetilde{\Phi}_{n, 1}^2,\dots, \widetilde{\Phi}_{n, m}^2)=
\int_{B(R)} G(|x|, |\Phi_{0, 1}|^2,\dots,|\Phi_{0, m}|^2).
\end{align}
Indeed, $\widetilde{\Phi}_{n, j}(x)    \rightarrow  \Phi_{0, j}(x) $ for a.e. $x \in  B(R)$,
and there exist $m$ functions
$f_j \in    L^{l_j +2}(B(R))$ such that $\widetilde{\Phi}_{n, j}(x)  \leq f_j(x)$ for a.e.
$x \in  B(R)$.  Of course $G(|x|, \widetilde{\Phi}_{n, 1}^2(x),\dots, \widetilde{\Phi}_{n, m}^2(x)) $
converges pointwise to $G(|x|, |\Phi_{0, 1}|^2(x),\dots,|\Phi_{0, m}|^2(x))$ in $  B(R)$ and, from ($G1$),
\begin{align*}
 G(|x|, \widetilde{\Phi}_{n, 1}^2,\dots, \widetilde{\Phi}_{n, m}^2)\leq K \Big(\sum_{j=1}^m f_j^2 +\sum_{j=1}^m f_j^{l_j+2} \Big)\in L^1(B(R)),
\end{align*}
Assertion~\eqref{GBR} then simply follows by dominated convergence.
Fixed $\eps > 0$, in light of~\cite[Lemma A.IV]{BerestyckiLions} and assumption ($G2$),  there
exist  $R(\eps) \geq R_0 > 0$ and $S_0 > 0$ such that, for  all $|x| > R(\eps)$,
$\widetilde{\Phi}_{n, j}(x) <  S_0$ for every $j=1, \dots, m$ and for all $n \in \N$.
Hence, by ($G2$), we have
\begin{align*}
  \int _{ B^c(R(\eps))} G(|x|, \widetilde{\Phi}_{n, 1}^2,\dots, \widetilde{\Phi}_{n, m}^2) \leq
  \eps \sum_{j=1}^m \int_{ B^c(R(\eps))}  \widetilde{\Phi}_{n, j}^2(x) \leq \eps c.
 \end{align*}
Now, observe that,  since $\widetilde{\Phi}_{n, j} (x)    \to \Phi_{0, j}(x) $ a.e.,
also $\widetilde{\Phi}_{0, j}(x) <  S_0$  for all $|x| > R(\eps)$. Then recalling that
also $\int \widetilde{\Phi}_{0, j}^2 \leq c$, we obtain
$$
  \int _{ B^c(R(\eps))} G(|x|, \widetilde{\Phi}_{0, 1}^2,\dots, \widetilde{\Phi}_{0, m}^2) \leq  \eps c.
$$
By means of~\eqref{GBR},  choosing $R=R(\eps)$, there exists
 $\nu_{\eps} \in \N$ such that, for all  $n \geq \nu_{\eps} $
\begin{align*}
\left| \int_{B(R(\eps))} G(|x|, \widetilde{\Phi}_{n, 1}^2,\dots, \widetilde{\Phi}_{n, m}^2)-
\int_{B(R(\eps))}  G(|x|, \Phi_{0, 1}^2,\dots,\Phi_{0, m}^2) \right| < \eps.
\end{align*}
Hence~\eqref{convergG} is proved too. Finally, we come to the proof of~\eqref{convergW}. We know that,
since $\widetilde{\Phi}_{n,j} $  is a sequence of radial functions, bounded in $H^1(\R^N)$,  by~\cite[Theorem A.I']{BerestyckiLions},
up to a subsequence, $\widetilde{\Phi}_{n, j} \to  \Phi_{0, j} $ strongly  in $L^{\hat q \mu}(\R^N)$ as $n\to\infty$,
where $\hat q= \frac{2q}{2q-1}$ and
$2 < \hat q \mu < 2^*$. Then, there exists a function $a_j \in L^{\hat q \mu}(\R^N) $  such that,
$ \widetilde{\Phi}_{n, j}(x)  \leq a_j(x)$ for a.e.~$x \in \R^N$. By the continuity of $h$, for
a.e.~$x, y \in \R^N$ we have
\begin{align*}
\lim_{n \rightarrow \infty}  W_{ij}(|x-y|)  h(|\widetilde{\Phi}_{n, i}(x)|)  h(|\widetilde{\Phi}_{n, j}(y)|)=
     W_{ij}(|x-y|)  h(|\Phi_{0, i}(x)|)  h(|\Phi_{0, j}(y)|).
\end{align*}
Furthermore, since $h$ is non-decreasing, we have for a.e. $x, y \in \R^N$
\begin{align*}
W_{ij}(|x-y|)  h(|\widetilde{\Phi}_{n, i}(x)|)  h(|\widetilde{\Phi}_{n, j}(y)|) \leq
W_{ij}(|x-y|)  h(a_i(x))  h(a_ j(y))
\end{align*}
where the right hand side function is in $L^1(\R^{2N})$ by means of Hardy-Littlewood Sobolev inequality
\begin{align*}
\iint W_{ij}(|x-y|)  h(a_i(x))  h(a_ j(y)) dx dy \leq \| a_i \|_{L^{\hat q \mu }(\R^N)}^{\mu} \| W_{ij} \|_{L^q_w(\R^N) }\| a_j \|_{L^{\hat q \mu }(\R^N)}^{\mu}.
\end{align*}
Then, by \eqref{gradSemicontInf}-\eqref{convergW}, \eqref{semicontFunct} is proved.
This yields $\mathcal E (\Phi_0) \leq I_c $, concluding the proof.
\end{proof}

\begin{proposition}
Assume conditions ($G3$), ($h0$) and ($h1$). If   $I_c < 0$, then    $\mathcal E (\Phi_0) = I_c$   for every $c > 0$.
\end{proposition}

\begin{proof}
In view of Proposition \ref{SchwarzLimit}, we know that $\mathcal E (\Phi_0) \leq  I_c$
and $\| \Phi_0 \|^2_{\mathcal{L}^2} \leq c$.
It is sufficient to prove that $\Phi_0 \in {\mathcal S}_c$.
First, we observe that,  by ($G$) and ($h$),  $\mathcal E (\textbf{0}) =0$ then $\Phi_0 \neq \textbf{0}$.
 Otherwise, by the negativity assumption on $I_c$, we would have
$$
0 = \mathcal E (\Phi_0) \leq I_c < 0,
$$
then a contradiction. Define $t = \frac{c^{1/2}}{\| \Phi_0 \|_{\mathcal{L}^2}}$,
we have that $t \Phi_0 \in {\mathcal S}_c$ and, by
$\| \Phi_0 \|^2_{\mathcal{L}^2} \leq c$,  $t \geq 1$. So, by ($G3$), ($h1$) and Proposition \ref{SchwarzLimit},
we have that
\begin{align*}
 \mathcal E (t \Phi_0) &  =  \frac{1}{2} \sum_{j=1}^ m \| \nabla (t \Phi_{0, j})  \|^2_{L^2}
- \frac{1}{2} \sum_{j=1}^ m V(x) \| t \Phi_{0, j}\|^2_{L^2}- \int G(|x|, t^2 \Phi_{0, 1}^2,\dots, t^2 \Phi_{0, m}^2)  dx   \\ \nonumber
& -  \frac12 \sum_{i, j=1}^{m} \iint W_{ij}(|x-y|)  h( t \Phi_{0, i}(x))  h(t \Phi_{j}(y))  dx dy  \leq t^2 \mathcal E ( \Phi_0) \leq t^2 I_c.
\end{align*}
Thus, $I_c \leq t^2 I_c$ and, by the negativity assumption on $I_c$, we have that $t \leq 1$. Hence, $t=1$ and by the definition of $t$,   $\| \Phi_0 \|^2_{\mathcal{L}^2} =  c$ thus proving the thesis.

\end{proof}

\section{Negativity of $I_c$}

The following results provides sufficient conditions in order to get the condition
that the minimum value is negative for all values of $c$.

\begin{proposition}
Assume conditions ($V0$), ($W1$) and either condition ($G4$) or condition
($h2$). Then $I_c < 0$ for all $c > 0$.
\end{proposition}

\begin{proof}
    In the following we shall assume both ($G4$) and ($h2$). It will be clear
    by the argument that follows that only one of these assumptions is actually
    sufficient to provide the desired conclusion.
Given $c>0$, we fix a positive function $\phi$ in $L^{\infty}(\R^N)$
such that $\|\phi\|_{L^2}^2= c$. Then, setting $\Phi = (\phi, 0, \dots, 0) \in \mathcal{H}^1$,
of course we, have $\Phi \in {\mathcal S}_c$. Now, for all $0 < t < 1$, let us
define $\phi_t(x)= t^{N/2} \phi(tx)$ and set $\Phi_t(x)=  (\phi_t(x), 0, \dots, 0)$.
Clearly, $\| \phi_t  \|_{L^2}^2= c$ and $\Phi_t \in {\mathcal S}_c$, for all $0<t<1$.
If we now evaluate the energy functional ${\mathcal E}$ at $\Phi_t$,
by a  change of variable and exploiting the assumptions,
for every $0<t< \min \{ t_1,  \frac{1}{R_2}\}$  sufficiently small, we have that
$$
0\leq t^{N/2} \phi(x)\leq t^{N/2} \|\phi\|_{L^\infty}\leq S_1,
\quad
0\leq t^N \phi^2(x)\leq t^N \|\phi\|^2_{L^\infty}\leq S_2,
$$
with $S_1$, $S_2$ and $R_2$ in assumptions  ($G4$) and ($h2$) so that
\begin{align*}
\mathcal E( \Phi_t) & =\frac12  \int \left| \nabla \phi_t (x) \right|^2  dx
- \frac12 \int V(|x|) \phi_t^2(x) \, dx- \int G(|x|, \phi_t^2(x), 0, \dots, 0)  dx   \\
& -  \frac12  \iint W_{11}(|x-y|)  h(\phi_{t}(x))  h(\phi_{t}(y))  dx dy \\
& = \frac{t^2}{2} \int \left| \nabla \phi (x) \right|^2  dx  - \frac12 \int V \left(\frac{|x|}{t} \right) \phi^2(x) \, dx
- t^{-N}\int G\Big(\frac{|x|}{t}, t^N \phi^2(x), 0, \dots, 0 \Big)  dx   \\ \nonumber
& - \frac{t^{-2N}}{2}  \iint W_{11} \left(\frac{|x-y|}{t} \right)  h( t^{N/2} \phi(x))  h( t^{N/2} \phi(y)) \,  dx dy  \\
& \leq \frac{t^2}{2} \int \left| \nabla \phi (x) \right|^2  dx
- t^{-N}\int_{\{|x|\geq  1 \}} G\Big(\frac{|x|}{t}, t^N \phi^2(x), 0, \dots, 0 \Big)  dx   \\ \nonumber
& - \frac{t^{-2N}}{2}  \iint W_{11} \left(\frac{|x-y|}{t} \right)  h( t^{N/2} \phi(x))  h( t^{N/2} \phi(y)) \,  dx dy  \\
& \leq Dt^2- E t^{-N} t ^{N \gamma}
- F t^{-2N} t^{\Gamma} t^{N \beta} ,
\end{align*}
where we have set
$$
D:=\frac{1}{2} \|\nabla \phi  \|_{L^2}^2,\quad
E:= B \int_{\{|x| \geq 1 \}} \phi^{2 \gamma}dx ,\quad
F:= A^2 C \iint \frac{\phi^{\beta}(x) \phi^{\beta}(y)}{|x-y|^{\gamma}} dx dy.
$$
In conclusion, for $t$ small enough, we get
\begin{equation*}
I_c\leq {\mathcal E}( \Phi_t) \leq   t^2 \big( D- E t ^{N \gamma-N-2}- F  t^{\Gamma +N \beta -2N-2} \big),
\end{equation*}
where, by the assumptions of $\gamma,\beta$ and $\Gamma$,
$$
N \gamma-N-2 < 0 \quad \hbox{and} \quad   \Gamma +N \beta -2N-2 < 0.
$$
By taking $t>0$ sufficiently small,  we have that $I_c \leq \mathcal E( \Phi_t)  < 0$, proving the assertion.
\end{proof}

\begin{remark}
    Notice that, if $W$ is a typical convolution kernel of the form $W(x)=|x|^{-\Gamma}$, it follows
    that $W$ belongs to the space $L^{q}_w(\R^N)$ where $q=\frac{N}{\Gamma}$. Moreover, thinking about the
    important model situation $h(s)=s^\mu$, we have $\beta=\mu$. Then, we have
    $$
    \Gamma +N \beta -2N-2 < 0
    \quad
    \Leftrightarrow
    \quad
    \frac{N}{q} +N \mu -2N-2 < 0
    \quad
    \Leftrightarrow
    \quad
    \mu<2-\frac{1}{q}+\frac{2}{N},
    $$
    which is the condition on $h$ we are already familiar with.
\end{remark}

\bigskip
\bigskip

\end{document}